\documentclass[twocolumn]{autart}    

\usepackage{amsmath,amssymb,epsf,epsfig,times}
\usepackage[all]{xy}
\usepackage{graphicx,color}
\usepackage{subfigure}
\usepackage{url}
\usepackage{cite}
\usepackage{epstopdf}
\usepackage{epsfig} 

\newtheorem{theorem}{Theorem}[section]
\newtheorem{lemma}{Lemma}[section]

\def\QED{\mbox{\rule[0pt]{1.5ex}{1.5ex}}}
\def\endproof{\hspace*{\fill}~\QED\par\endtrivlist\unskip}
\newcommand{\re}{\mathbb{R}}

\newtheorem{assumption}[theorem]{Assumption}

\newtheorem{remark}[theorem]{Remark}

\newcommand{\OMIT}[1]{}
\allowdisplaybreaks[4]

\begin{document}
	\begin{frontmatter}

		\title{Distributed Average Tracking for Double-integrator Multi-agent Systems with Reduced Requirement on Velocity Measurements}
		
		
		\author[Paestum]{Sheida Ghapani}\ead{sghap001@ucr.edu}~,
		\author[Paestum]{Wei Ren}\ead{ren@ee.ucr}~,
		\author[Paestum1]{Fei Chen}\ead{feichen@xmu.edu.cn}~,  
		\author[Paestum2]{Yongduan Song}\ead{ydsong@cqu.edu.cn}

		\address[Paestum]{Department of Electrical and Computer Engineering, University of California, Riverside, CA 92521, USA}
		\address[Paestum1]{Department of Automation, Xiamen University, Xiamen, China, 361005}
		\address[Paestum2]{School of Automation, Chongqing University, Chongqing, China, 400044}
		


		\begin{abstract}
			
			This paper addresses distributed average tracking for a group of physical double-integrator agents under an undirected graph with reduced requirement on velocity measurements.
			The idea is that multiple agents track the average of multiple time-varying input signals, each of which is available to only one agent, under local interaction with neighbors.
			We consider two cases.
			First, a distributed discontinuous algorithm and filter are proposed
			, where each agent needs the relative positions between itself and its neighbors and its neighbors' filter outputs obtained through communication but the requirement for either absolute or relative velocity measurements is removed.
			The agents' positions and velocities must be initialized correctly, but the algorithm can deal with a wide class of input signals with bounded acceleration deviations.
			Second, a distributed discontinuous algorithm and filter are proposed to remove the requirement for communication and accurate initialization.
			Here each agent needs to measure the relative position between itself and its neighbors and its own velocity but the requirement for relative velocity measurements between itself and its neighbors is removed.
			The algorithm can deal with the case where the input signals and their velocities and accelerations are all bounded. Numerical simulations are also presented to illustrate the theoretical results.

		\end{abstract}
		
	\end{frontmatter}
	\section{INTRODUCTION}
	
	This paper studies the following distributed average tracking problem: given a group of agents and one time-varying input signal per each agent, design a control law for the agents based on local information such that all the agents will finally track the average of these input signals.
	The problem has found applications in distributed sensor fusion \cite{olfati2004consensus}, feature-based map merging \cite{aragues2012distributed}, distributed Kalman filtering \cite{bai2011}, where distributed computation of multiple time-varying signals are involved.
	Compared with the consensus problem, distributed average tracking poses more theoretical challenges, since the tracking objective is time-varying and is not available to any agent.
	
	In the literature, linear distributed algorithms have been employed for special kinds of time-varying input signals. 
	Ref. \cite{spanos2005dynamic} uses frequency domain analysis to study consensus on the average of multiple input signals with
	steady-state values.
	In \cite{freeman2006}, a proportional algorithm and a proportional-integral algorithm are proposed to achieve distributed average tracking with a bounded tracking
	error, where accurate estimator initialization is relaxed in the proportional-integral algorithm.
	In \cite{bai2010}, the internal model principle is employed to extend the proportional-integral algorithm to a special group of time-varying input signals with a common denominator in their Laplace transforms, where the denominator also needs to be used in the estimator  design.
	Ref. \cite{Zhumartinez2010} addresses discrete-time distributed average tracking of time-varying input signals whose $n$th order difference is bounded with a bounded error.
	In \cite{kiaDACsingularity}, the authors propose 1st-order-input and 2nd-order-input consensus algorithms to allow the agents to track the average of their dynamic
	input signals with a pre-specified rate, where the interaction is described by a strongly connected and weight-balanced directed graph.
	In \cite{kiaauthority}, a distributed continuous average tracking problem is addressed with some steady-state errors in , where the privacy of each agent's input signal is preserved.
	
	However, linear algorithms cannot ensure distributed average tracking for general input signals. Therefore, some researchers employ nonlinear tracking algorithms.
	In \cite{Nosrati2012}, a class of nonlinear algorithms is proposed for input signals with bounded deviations, where the tracking error is proved to be bounded.
	A nonsmooth algorithm is proposed in \cite{chen2012distributed}, which is able to track arbitrary time-varying input signals with bounded derivatives.
	All the above references primarily study the distributed average tracking problem from a distributed estimation perspective, where the agents implement local estimators through communication with neighbors freely without the need for obeying certain physical agent dynamics.
	However, there are applications where the distributed average tracking problem is relevant for designing distributed control laws for physical agents.
	One application is the region-following formation control \cite{cheah2009region}, where a swarm of robots is required to move inside a dynamic region while keeping a desired formation.
	Then the dynamics of the physical agents must be taken into account in the control law design and the dynamics themselves introduce further challenges to the distributed average tracking problem.
	For example, the control law designed for physical agents with single-integrator dynamics can no longer be used directly for physical agents subject to more complicated dynamic equations.
	Distributed average tracking for physical agents with double-integrator dynamics is studied in \cite{feidoubleintegrator}, where the input signals are allowed to have bounded accelerations.
	Distributed average tracking for physical agents with general linear dynamics is addressed in \cite{DBLP:journals/corr/ZhaoDL13} and \cite{FeiRobust}.
	Ref. \cite{FeiRobust} proposes a discontinuous algorithm, while a continuous algorithm is employed in \cite{DBLP:journals/corr/ZhaoDL13} with, respectively, static and adaptive coupling strengths.
	Also in \cite{fei2015EulerDAC} a proportional-integral control scheme is extended to achieve distributed average tracking for physical Euler-Lagrange systems for two different kinds of input signals with steady states and with bounded derivatives.
	
	It is noted that in \cite{feidoubleintegrator} both relative position and relative velocity measurements are required in the control laws designed for double-integrator agents.
	While double-integrator dynamics can be viewed as a special case of general linear dynamics, the distributed average tracking algorithms in \cite{FeiRobust,fei2015EulerDAC}, when applied to double-integrator systems, still need relative position and velocity measurements.
	In practice, velocity measurements are usually less accurate and more expensive than position measurements.
	In addition, relative velocity measurements are often more challenging and expensive than absolute velocity measurements.
	We are hence motivated to solve the distributed average tracking problem for physical double-integrator agents with reduced requirement on velocity measurements. This paper expands on our preliminary work reported in \cite{GhapaniRenChen15_ACC}.
	In the context of distributed average tracking, reducing velocity measurements poses significant theoretical challenges. The reason is that unlike the consensus or single-leader coordinated tracking problems, there are significant additional inherent challenges in distributed average tracking as none of the agents has the tracking objective available.
	
	The contribution of this paper is summarized as follows. Two distributed algorithms (controller design combined with filter design) are introduced to achieve distributed average tracking with reduced requirement on velocity measurements.
	Each algorithm has its own relative benefits and is feasible for different application scenarios.
	In the first algorithm design, there is no need for either absolute or relative velocity measurements.
	Each agent's algorithm only employs its local relative positions with respect to neighbors, its neighbors' filter outputs accessed through communication, and the acceleration of its own input signal.
	The algorithm allows the agents to track the average of a large class of time-varying input signals with bounded acceleration deviations, provided that the agents are correctly initialized.
	In the second algorithm design, there is no need for correct initialization and relative velocity measurements. Also inter-agent communication is not necessary and the algorithm can be implemented using only local sensing, which is desirable for certain applications (e.g., deep-space spacecraft formation flying) where communication might not be desirable or available.
	Each agent's algorithm only employs its local relative positions with respect to neighbors, its own velocity, and its own input signal.
	Distributed average tracking can be achieved provided that the input signals and their velocities and accelerations are all bounded.
	
	{\it Notations:} Throughout the paper, $\mathbb{R}$ denotes the set of all real numbers and $\mathbb{R}^+$ the set of all positive real numbers. Let $\mathbf{1}_n$
	and $\mathbf{0}_n$ denote the $n \times 1$ column vector of all ones and all zeros respectively.
	Let $\lambda_{\max} (\cdot)$ and $\lambda_{\min} (\cdot)$ denote, respectively, the maximum and minimum eigenvalues of a square real matrix with real eigenvalues. %
	We use $\otimes$ to denote the Kronecker product, and $\mbox{sgn}(\cdot)$ to denote the $\mbox{signum}$ function defined componentwise.
	For a vector function ${x(t):\re\mapsto\re^m}$, define
	$||x(t)||_p$ as the p-norm, ${x(t)\in\mathbb{L}_2}$ if
	$\int_{0}^{\infty}x(\tau)^T x(\tau)\mbox{d}\tau<\infty$ and
	${x(t)\in\mathbb{L}_{\infty}}$ if for each element of $x(t)$, ${\text{sup}_{t \geq 0}|x_i(t)|<\infty}$, $i=1,\ldots,m$.

	\section{Problem Statement}
	Here we consider $n$ physical agents described by double-integrator dynamics
	\begin{align} \label{non-agents-dynamic}
		\dot{x}_i(t)=&v_i(t), \notag  \\
		\dot{v}_i(t)=&u_i(t),  \qquad i=1,\ldots,n,
	\end{align}
	where $x_i(t) \in \mathbb{R}^p$ and $v_i(t) \in \mathbb{R}^p$ are, respectively, agent $i$'s position and velocity, and $u_i(t)$ is its control input. Let ${x(t)=[x_1^T(t),\ldots,x_n^T(t) ]^T}$ and ${v(t)=[v_1^T(t),\ldots,v_n^T(t) ]^T}$.
	
	An \textit{undirected} graph $G \triangleq (V,E)$ is used to characterize the interaction topology among the agents, where ${V \triangleq \{1,\ldots,n\}}$ is the node set and $E \subseteq V \times V$ is the edge set.
	An edge $(j,i) \in E$ means that node $i$ can obtain information from node $j$ and vice versa.
	Self edges $(i,i)$ are not considered here.
	Let $m$ denote the number of edges in $E$, where the edges $(j,i)$ and $(i,j)$ are counted only once.
	The \textit{adjacency matrix} ${\mathbf{A}=[a_{ij}] \in \mathbb{R}^{n \times n}}$ of the graph $G$ is defined such that the edge weight ${a_{ij}=1}$ if ${(j,i) \in E}$ and ${a_{ij}=0}$ otherwise. For an undirected graph, ${a_{ij}=a_{ji}}$.
	The \textit{Laplacian matrix} ${L=[l_{ij}] \in \mathbb{R}^{n \times n}}$ associated with $\mathbf{A}$ is defined as ${l_{ii}=\sum_{j \ne i} a_{ij}}$ and ${l_{ij}=-a_{ij}}$, where ${i \ne j}$.
	For an undirected graph, $L$ is symmetric positive semi-definite.
	By arbitrarily assigning an orientation for the edges in $G$, let $D \triangleq [d_{ij}] \in  \mathbb{R}^{n \times m}$ be the \textit{incidence matrix} associated with $G$, where $d_{ij} = -1$ if the edge $e_j$ leaves node $i$, $d_{ij} = 1$ if it enters node $i$, and $d_{ij} = 0$ otherwise.
	The Laplacian matrix $L$ is then given by $L=DD^T$ \cite{GodsilRoyle01}.

	\begin{assumption} \label{conn-graph}
		The undirected graph $G$ is connected.
	\end{assumption}
	
	\begin{lemma} \cite{GodsilRoyle01} \label{eigen}
		Under Assumption \ref{conn-graph}, the Laplacian matrix $L$ has a simple zero eigenvalue such that $0=\lambda_1(L)<\lambda_2(L) \leq \ldots \leq \lambda_n(L)$, where $\lambda_i(\cdot)$ denotes the $i$th eigenvalue. Furthermore, for any vector $y \in \mathbb{R}^n$ satisfying ${\mathbf{1}_n^T y=0}$, $\lambda_2(L) y^Ty \leq y^T L y \leq \lambda_n(L) y^Ty$.
	\end{lemma}
	
	Suppose that each agent has a time-varying input signal $r_i(t) \in \mathbb{R}^p$, $i=1,\ldots,n$, satisfying
	\begin{align} \label{ref-dynamic}
		\dot{r}_i(t)=& v_i^r(t),   \notag \\
		\dot{v}_i^r(t)=& a_i^r(t),
	\end{align}
	where $v_i^r(t)$ and $a_i^r(t) \in \mathbb{R}^p$ are, respectively, the agent $i$'s input velocity and input acceleration.
	Define ${r(t)=[r_1^T,\ldots,r_n^T ]^T}$, ${v^r(t)=[{v_1^r}^T,\ldots,{v_n^r}^T ]^T}$ and ${a^r(t)=[{a_1^r}^T,\ldots,{a_n^r}^T ]^T}$. Here we assume that the input signals are generated internally by the agents and each agent has access to its own input signal, input velocity, and input acceleration.

	We study the distributed average tracking problem for the double-integrator agents. The goal is to design $u_i(t)$ for each agent to track the average of the input signals and input velocities, i.e.,
	\begin{align}
		\lim \limits_{t \to \infty} ||x_i(t)-\frac{1}{n} \sum_{j=1}^n r_j(t)||_2=&0, \notag \\
		\lim \limits_{t \to \infty} ||v_i(t)-\frac{1}{n} \sum_{j=1}^n v_j^r(t)||_2=&0, \qquad i=1,\ldots,n, \notag
	\end{align}
	where each agent has access to its own input information and has only local interaction with its neighbors. We are interested in controller design with reduced requirement on velocity measurements. First, we consider the scenario that each agent has communication capabilities but without the need for either absolute or relative velocity measurements. The agents achieves distributed average tracking for a wide class of input signals with bounded deviations among input accelerations in the presence of correct state initialization.
	Second, we consider the scenario that each agent has sensing but not necessarily communication capabilities without the need for relative velocity measurements. The agents achieve distributed average tracking even in the absence of correct state initialization.
	
	\section{Velocity free distributed average tracking in the presence of correct initialization}\label{corr-ini}
	In this section, we consider the case where the agents can initialize their initial conditions correctly. To remove both absolute and relative velocity measurements, we introduce the following filter for each agent
	\begin{align} \label{vel-free-fil-in}
		\dot{\upsilon}_i(t)=&\sum_{j \in N_i} [x_i(t)-x_j(t) ]-\beta \sum_{j \in N_i} [w_i(t)-w_j(t) ] \notag \\
		&-\gamma \sum_{j \in N_i} \mbox{sgn}[w_i(t)-w_j(t)]-a_i^r(t),
	\end{align}
	\begin{align} \label{vel-free-fil}
		w_i(t)=\upsilon_i(t) -\alpha \sum_{j \in N_i} [x_i(t)-x_j(t)], \qquad i=1,\ldots,n,
	\end{align}
	where $\upsilon_i \in \mathbb{R}^p$ is an auxiliary filter variable, $w_i \in \mathbb{R}^p$ is the filter output, and $\beta$, $\alpha$ and $\gamma \in \mathbb{R}^+$ are control gains to be designed.

	We propose the following distributed control law for agent $i$
	\begin{align} \label{vel-free-control}
		u_i(t)=&-\alpha \sum_{j \in N_i} [x_i(t)-x_j(t) ]+\beta \sum_{j \in N_i} [w_i(t)-w_j(t) ] \notag \\
		& +\gamma \sum_{j \in N_i} \mbox{sgn}[w_i(t)-w_j(t)] \notag \\
		&+a_i^r(t),  \qquad i=1,\ldots,n.
	\end{align}
	For notational simplicity, we will remove the index $t$ from variables in the reminder of the paper.
	
	\begin{assumption} \label{bounded-a}
		The deviations among input accelerations are bounded, i.e., ${\text{sup}_{t \in [0,\infty)} ||a_i^r(t)-a_j^r(t)||_2 \leq \bar{a}^r_d}$, $i$, ${j=1,\ldots,n}$, where $\bar{a}^r_d \in \mathbb{R}^+$.
	\end{assumption}
	
	\begin{assumption} \label{initial}
		\begin{align}
			\sum_{i=1}^n x_i(0)=&\sum_{i=1}^n r_i(0), \notag \\
			\sum_{i=1}^n v_i(0)=&\sum_{i=1}^n v_i^r(0). \notag \footnotemark
		\end{align}
		\footnotetext{A special choice is $x_i(0)=r_i(0)$, $v_i(0)=v_i^r(0)$ for all ${i=1,\ldots,n}$.}
	\end{assumption}

	\begin{theorem} \label{thm:DAC-vel-dynamic}
		Using the control law given by \eqref{vel-free-control}, \eqref{vel-free-fil}, and \eqref{vel-free-fil-in} for  system \eqref{non-agents-dynamic}, distributed average tracking is achieved asymptotically, provided that Assumptions \ref{conn-graph}, \ref{bounded-a} and \ref{initial} hold and the gains $\alpha$, $\gamma$ and $\beta$ are chosen such that
		\begin{align}
			\alpha > \max \{ 1, \frac{1}{\lambda_2(L)}, \frac{ \lambda_2(L)+1}{2 \lambda_2(L)} \}, \notag
		\end{align}
		\begin{align}
			\gamma > (n-1)\bar{a}^r_d,\notag
		\end{align}
		\begin{align}
			\beta > \frac{1+\alpha^4  \lambda_n^2(L) }{4 \lambda_2(L)(\alpha \lambda_2(L)-1)(\alpha-1)},\notag
		\end{align}
		where $\lambda_2(L)$ and $\lambda_n(L)$ are defined in Lemma \ref{eigen}.
	\end{theorem}
	
	\emph{Proof}:
	The proof contains two steps. First, we prove that for each agent, ${x_i-\frac{1}{n} \sum_{j=1}^n x_j \to 0}$ and ${v_i-\frac{1}{n} \sum_{j=1}^n v_j \to 0}$.
	Then by showing that $\sum_{i=1}^n x_i-\sum_{i=1}^n r_i \to 0$ and $\sum_{i=1}^n v_i-\sum_{i=1}^n v_i^r \to 0$, it can be concluded that ${x_i-\frac{1}{n} \sum_{j=1}^n r_j \to 0}$ and ${v_i-\frac{1}{n} \sum_{j=1}^n v_j^r \to 0}$, ${i=1,\ldots,n}$, and hence distributed average tracking is achieved.
	
	Using the control law \eqref{vel-free-control} for \eqref{non-agents-dynamic}, we can get
	\begin{align} \label{closed-loop}
		\dot{x}_i=&v_i, \notag  \\
		\dot{v}_i=&-\alpha \sum_{j \in N_i}[x_i(t)-x_j(t) ]+\beta \sum_{j \in N_i} [w_i(t)-w_j(t) ] \notag \\
		& +\gamma \sum_{j \in N_i} \mbox{sgn}[w_i(t)-w_j(t)] +a_i^r(t),  \qquad i=1,\ldots,n.
	\end{align}
	Define ${\xi_x=[M \otimes I_p] x}$, ${\xi_v=[M \otimes I_p] v}$ and ${\xi_w=[M \otimes I_p] w}$, where ${M=I_n-\frac{1}{n} \mathbf{1}_n \mathbf{1}_n^T} $ and $w=[w_1^T,\ldots,w_n^T]^T$.
	Then the closed-loop dynamics \eqref{closed-loop} can be rewritten as
	\begin{align} \label{cons-error}
		\dot{\xi}_x =&\xi_v, \notag \\
		\dot{\xi}_v=&-\alpha [L \otimes I_p] \xi_x+\beta [L \otimes I_p]\xi_w \notag \\
		&+\gamma [D \otimes I_p] \mbox{sgn}( [D^T \otimes I_p]\xi_w) +[M \otimes I_p] a^r,
	\end{align}
	where it can be easily proved that $M \times M=M$, $L \times M=L$ and $D \times M=D$.
	It follows from \eqref{vel-free-fil-in} and \eqref{vel-free-fil} that
	\begin{align} \label{vel-filter-compact}
		\dot{\xi}_w=& [L \otimes I_p] \xi_x-\beta [L \otimes I_p]\xi_w \notag \\
		&-\gamma [D \otimes I_p] \mbox{sgn}( [D^T \otimes I_p]\xi_w)-a^r- \alpha [L \otimes I_p] \xi_v.
	\end{align}
	Consider the following function
	\begin{align} 
		V=&\frac{1}{2}
		\xi^T
		\begin{bmatrix}
			\mu  [L \otimes I_p] &  I_{np} & I_{np} \\
			I_{np} & I_{np} &  I_{np} \\
			I_{np} & I_{np} &  \alpha I_{np}
		\end{bmatrix}
		\xi, \notag
	\end{align}
	where $\mu \in \mathbb{R}^+$ and ${\xi=
		\begin{bmatrix}
		\xi_x^T & \xi_v^T & \xi_w ^T
		\end{bmatrix}}^T$. Since $\xi_x^T \mathbf{1}_{np}=0$, using Lemma \ref{eigen}, we have
	\begin{align}
		V \geq &\frac{1}{2}
		\xi^T
		(\begin{bmatrix}
			\mu  \lambda_2(L)  &  1 & 1 \\
			1 & 1 &  1 \\
			1 & 1 &  \alpha
		\end{bmatrix}
		\otimes I_{np})
		\xi \notag
	\end{align}
	It can be proved that if $\mu > \frac{1}{\lambda_2(L)}$ and $\alpha>1$, then
	$\begin{bmatrix}
	\mu  \lambda_2(L)  &  1 \\
	1 & 1
	\end{bmatrix} -\frac{1}{\alpha} \begin{bmatrix}
	1 \\
	1
	\end{bmatrix} \begin{bmatrix}
	1  &  1
	\end{bmatrix}>0$. Thus, using Shur complement it is concluded that $V$ is positive definite.
	The derivative of $V$ along the trajectories of \eqref{cons-error} and \eqref{vel-filter-compact} is given as
	\begin{align} \label{vel-lyp-der}
		\dot{V}=&
		\xi^T
		\begin{bmatrix}
			\mu  [L \otimes I_p] &  I_{np} & I_{np} \\
			I_{np} & I_{np} &  I_{np} \\
			I_{np} & I_{np} &  \alpha I_{np}
		\end{bmatrix}
		\dot{\xi} \notag \\
		=& \mu  \xi_x^T[L \otimes I_p] \xi_v+  \xi_v^T \xi_v+  \xi_w^T \xi_v- \alpha \xi_x^T [L \otimes I_p] \xi_x \notag \\
		&- \alpha \xi_v^T [L \otimes I_p] \xi_x+ \beta \xi_w^T [L \otimes I_p] \xi_w \notag \\
		&+\gamma \xi_w^T  [D \otimes I_p] \mbox{sgn}( [D^T \otimes I_p]\xi_w) +\xi_w^T a^r + \xi_x^T [L \otimes I_p] \xi_x  \notag \\
		&-\alpha \xi_x^T [L \otimes I_p] \xi_v+\xi_v^T [L \otimes I_p] \xi_x -\alpha \xi_v^T [L \otimes I_p] \xi_v\notag \\
		&-\alpha \beta  \xi_w ^T [L \otimes I_p] \xi_w -\alpha \gamma \xi_w ^T [D \otimes I_p] \mbox{sgn}( [D^T \otimes I_p]\xi_w) \notag \\
		&-\alpha \xi_w^T a^r-\alpha^2  \xi_w^T [L \otimes I_p] \xi_v,
	\end{align}
	where $[M \otimes I_p] \xi_x=[M^2 \otimes I_p]x=[M \otimes I_p]x=\xi_x$ and using the same analysis $[M \otimes I_p] \xi_v=\xi_v$ and $[M \otimes I_p] \xi_w=\xi_w$.
	We can also analyze the term $(1-\alpha) \xi_w^T a^r $ as
	\begin{align}\label{w-bound}
		(1-\alpha) \xi_w^T  a^r=& (1-\alpha) \xi_w^T [M^2 \otimes I_p] a^r  \notag \\
		\leq & (\alpha-1)||[M \otimes I_p]\xi_w||_2  ||[M \otimes I_p]a^r||_2 \notag \\
		\leq &  \frac{(\alpha-1) \bar{a}^r_d}{n} \sum_{i=1}^n \sum_{j=1,j \ne i}^n ||\xi_{wi}-\xi_{wj}||_2   \notag \\
		\leq & \frac{(\alpha-1) \bar{a}^r_d}{n} \sum_{i=1}^n \max\limits_{i}\{\sum_{j=1,j \ne i}^n||\xi_{wi}-\xi_{wj}||_2\}  \notag \\
		\leq &  (\alpha-1) \bar{a}^r_d \max\limits_{i}\{\sum_{j=1,j \ne i}^n ||\xi_{wi}-\xi_{wj}||_2\} \notag \\
		\leq &  \frac{(\alpha-1)(n-1)}{2}  \bar{a}^r_d  \sum_{i=1}^n \sum_{j \in N_i} ||\xi_{wi}-\xi_{wj}||_2  \notag \\
		\leq &  \frac{(\alpha-1)(n-1) }{2}\bar{a}^r_d  \sum_{i=1}^n \sum_{j \in N_i} ||\xi_{wi}-\xi_{wj}||_1,
	\end{align}
	where we have used the Assumption \ref{bounded-a} and $||\cdot||_2 \leq ||\cdot||_1$ to obtain second and last inequalities, respectively.
	
	Let ${\mu=2\alpha-1}$. Note that $\alpha>\frac{ \lambda_2(L)+1}{2 \lambda_2(L)}$ ensures that $\mu > \frac{1}{\lambda_2(L)}$. Combining \eqref{vel-lyp-der} and \eqref{w-bound}, we can get
	\begin{align}  \label{lyp-der-fi}
		\dot{V} \leq  &  \xi_v^T \xi_v+  \xi_w^T \xi_v+(1-\alpha) \xi_x^T [L \otimes I_p] \xi_x \notag \\
		&+(1-\alpha)\beta \xi_w^T [L \otimes I_p] \xi_w \notag \\
		&-\gamma (\alpha-1) \sum_{i=1}^n \sum_{j \in N_i} \xi_{wi}^T \mbox{sgn}(\xi_{wi}-\xi_{wj}) \notag \\
		&+\frac{(\alpha-1)(n-1) }{2}\bar{a}^r_d  \sum_{i=1}^n \sum_{j \in N_i} ||\xi_{wi}-\xi_{wj}||_1   \notag \\
		& -\alpha \xi_v^T [L \otimes I_p] \xi_v-\alpha^2  \xi_w^T [L \otimes I_p] \xi_v \notag \\
		=& \xi_v^T \xi_v+  \xi_w^T \xi_v+(1-\alpha) \xi_x^T [L \otimes I_p] \xi_x \notag \\
		&+(1-\alpha)\beta \xi_w^T [L \otimes I_p] \xi_w \notag \\
		&-\frac{\gamma (\alpha-1)}{2} \sum_{i=1}^n \sum_{j \in N_i} (\xi_{wi}-\xi_{wj})^T \mbox{sgn}(\xi_{wi}-\xi_{wj}) \notag \\
		&+\frac{(\alpha-1)(n-1) }{2}\bar{a}^r_d  \sum_{i=1}^n \sum_{j \in N_i} ||\xi_{wi}-\xi_{wj}||_1   \notag \\
		& -\alpha \xi_v^T [L \otimes I_p] \xi_v-\alpha^2  \xi_w^T [L \otimes I_p] \xi_v \notag \\
		\leq &\xi^T
		\begin{bmatrix}
			Q_{11} &  Q_{12} & Q_{13}  \\
			Q_{21} &  Q_{22} & Q_{23}  \\
			Q_{31} &  Q_{32} & Q_{33}
		\end{bmatrix}
		\xi,
	\end{align}
	where $Q_{11}=(1-\alpha) \lambda_2(L) I_{np}$, $Q_{22}=[1-\alpha \lambda_2(L)]I_{np}$, $Q_{23}=Q_{32}=I_{np}-\alpha^2 [L \otimes I_p]$, $Q_{33}=(1-\alpha) \beta \lambda_2(L) I_{np}$ $Q_{12}=Q_{13}=Q_{21}=Q_{31}= 0_{np}$,
	and we have used ${\gamma>(n-1)\bar{a}^r_d}$ and Lemma \ref{eigen} for the last inequality.
	We also used $\xi_w ^T [D \otimes I_p] \mbox{sgn}( [D^T \otimes I_p]\xi_w)=\sum_{i=1}^n \sum_{j \in N_i} \xi_{wi}^T \mbox{sgn}(\xi_{wi}-\xi_{wj})$ to derive the condition of $\gamma$.
	
	Therefore, if the control gains $\alpha$ and $\beta$ satisfy the constraint $\alpha > \max \{ 1, \frac{1}{\lambda_2(L)} \}$ and $\beta > \frac{1+\alpha^4  \lambda_n^2(L) }{4 \lambda_2(L)(\alpha \lambda_2(L)-1)(\alpha-1)}$, it can be seen that $\dot{V}$ is negative definite.
	Integrating both sides of \eqref{lyp-der-fi}, we can obtain $\xi_x$, $\xi_v$, $\xi_w \in \mathbb{L}_2$.
	Since $V$ is positive definite and $\dot{V}$̇ is negative definite, $\xi_x$, $\xi_v$, $\xi_w \in L_{\infty}$.
	Thus, we can get $\xi_x$, $\xi_v$, $\xi_w \in L_{2} \cap L_{\infty}$ and $\dot{\xi}_x$, $\dot{\xi}_v$, $\dot{\xi}_w \in L_{\infty}$. Using Barbalat's lemma \cite{slotine1991}, it is concluded that $\xi_x$, $\xi_v$, $\xi_w \to 0$ as $t \to \infty$ which imply that $x_i \to \frac{1}{n}\sum_{j=1}^n x_j$ and $v_i \to \frac{1}{n} \sum_{j=1}^n v_j$, $i=1,\ldots,n$.
	
	Second, we analyze the terms ${\sum_{i=1}^n x_i}$ and ${\sum_{i=1}^n v_i}$. Considering the closed-loop dynamics \eqref{closed-loop}, under Assumption \ref{conn-graph}, the derivative of ${\sum_{i=1}^n x_i}$ and ${\sum_{i=1}^n v_i}$ can be calculated as
	\begin{align}
		\sum_{i=1}^n \dot{x}_i=\sum_{i=1}^n v_i, \notag \\
		\sum_{i=1}^n \dot{v}_i=\sum_{i=1}^n a_i^r, \notag
	\end{align}
	where we have used the fact that the graph $G$ is undirected. Using Assumption \ref{initial}, it can be proved that ${\sum_{i=1}^n x_i = \sum_{i=1}^n r_i}$ and ${\sum_{i=1}^n v_i = \sum_{i=1}^n v_i^r}$.
	
	Combining the two parts shows that ${x_i \to \frac{1}{n}\sum_{j=1}^n r_j}$ and ${v_i \to \frac{1}{n}\sum_{j=1}^n v_j^r}$ asymptotically. Therefore, distributed average tracking is achieved asymptotically.
	\endproof

	\begin{remark} Compared with \cite{feidoubleintegrator}, using the algorithm defined by \eqref{vel-free-fil-in}-\eqref{vel-free-control}, the requirement for velocity measurements is removed. In addition, the only limitation on the input signals is that the deviations among their accelerations are bounded. Therefore, the allowable input signals are more general than those in \cite{feidoubleintegrator}.
	\end{remark}

	\section{Distributed average tracking in the absence of neighbors' velocity information and correct initialization} \label{abs-corr-ini}
	In the previous section, the proposed algorithm solves the distributed average tracking problem in the presence of communication and accurate state initialization without velocity measurements. In this section we deal with the distributed average tracking problem without communication in the absence of neighbors' velocity information and correct initialization.
	
	We introduce the following filter
	\begin{align} \label{ini-free-fil-in}
		\dot{\upsilon}_i=&\sum_{j \in N_i} a_{ij}(x_i-x_j) -\beta w_i  -\gamma \mbox{sgn}(w_i)- \kappa r_i \notag \\
		&-\kappa v_i^r-a_i^r,
	\end{align}
	\begin{align} \label{ini-free-fil}
		w_i=&\upsilon_i-\alpha \sum_{j \in N_i} (x_i-x_j) , \qquad i=1,\ldots,n,
	\end{align}
	where $\upsilon_i \in \mathbb{R}^p$ is an auxiliary filter variable, $w_i \in \mathbb{R}^p$ is the filter output,  and $\kappa$, $\alpha$, $\gamma$, $\beta \in \mathbb{R}^+$ are control gains.

	\begin{assumption} \label{bounded-v-a}
		Each input signal, its velocity and acceleration are bounded, i.e., ${\text{sup}_{t \in [0,\infty)} ||r_i(t)||_2 \leq \bar{r}}$, ${\text{sup}_{t \in [0,\infty)} ||v_i^r(t)||_2 \leq \bar{v}^r}$ and ${\text{sup}_{t \in [0,\infty)} ||a_i^r(t))||_2 \leq \bar{a}^r}$, $i=1,\ldots,n$, where $\bar{r}$, $\bar{v}^r$, $\bar{a}^r \in \mathbb{R}^+$.
	\end{assumption}
	
	We propose the following control input for each agent
	as
	\begin{align} \label{ini-free-control}
		u_i=&-\kappa (x_i-r_i )-\kappa (v_i-v_i^r )-\alpha \sum_{j \in N_i} (x_i-x_j) +\beta w_i \notag \\
		&+\gamma \mbox{sgn}(w_i) +a_i^r,  \qquad i=1,\ldots,n.
	\end{align}
	
	%
	
	The following theorem presents sufficient conditions to solve distributed average tracking without neighbors' velocity and accurate initialization.
	
	\begin{theorem} \label{thm:ini-free-DAC}
		Distributed average tracking is achieved for  system \eqref{non-agents-dynamic} using \eqref{ini-free-control}, \eqref{ini-free-fil}  and \eqref{ini-free-fil-in} asymptotically, under Assumptions \ref{conn-graph} and \ref{bounded-v-a} and provided that the gains $\kappa$, $\alpha$, $\gamma$ and $\beta$  are chosen such that
		\begin{align}
			\kappa > 1 , \notag
		\end{align}
		\begin{align}
			\alpha >  \sqrt{n}, \notag
		\end{align}
		\begin{align}
			\gamma > \frac{\alpha+1}{\alpha-\sqrt{n}}(\kappa\bar{r}+\kappa \bar{v}^r+\bar{a}^r),\notag
		\end{align}
		\begin{align}
			\beta >& \frac{1}{4(\alpha -1)}[\frac{\kappa^2}{\kappa +(\alpha -1) \lambda_2(L)} \notag \\
			&+\frac{(\kappa-1)^2 + 2 \alpha^2 (\kappa-1) \lambda_n(L) + \alpha^4 \lambda_n^2(L)}{\alpha \lambda_2(L) + \kappa -1}],\notag
		\end{align}
		where $\lambda_2(L)$ and $\lambda_n(L)$ are defined in Lemma \ref{eigen}.
	\end{theorem}
	
	\emph{Proof}:
	Similar to Section III, the proof contains two parts.
	In first part, we prove that ${x_i-\frac{1}{n} \sum_{j=1}^n x_j \to 0}$, ${v_i-\frac{1}{n} \sum_{j=1}^n v_j \to 0}$ and $ {w_i \to 0}$ asymptotically.
	Define $\xi_x$, $\xi_v$, $M$, and $w$ as in the proof of Theorem \ref{thm:DAC-vel-dynamic}.
	Now rewrite the dynamics \eqref{non-agents-dynamic} using \eqref{ini-free-control} in vector form as
	\begin{align} 
		\dot{\xi}_x =&\xi_v, \notag \\
		\dot{\xi}_v=& -\kappa \xi_x+\kappa [M \otimes I_p]r-\kappa \xi_v+\kappa [M \otimes I_p]v^r \notag \\
		&-\alpha [L \otimes I_p] \xi_x  +\beta [M \otimes I_p] w+ \gamma [M \otimes I_p] \mbox{sgn}(w) \notag \\
		&+[M \otimes I_p]a^r. \notag
	\end{align}
	Consider the following Lyapunov function candidate
	\begin{align} 
		V=&\frac{1}{2}
		\begin{bmatrix}
			\xi_x^T &  \xi_v^T & w^T
		\end{bmatrix}
		\begin{bmatrix}
			\mu_1  [L \otimes I_p] + \mu_2 I_{np} &  I_{np} & I_{np} \\
			I_{np} & I_{np} &  I_{np} \\
			I_{np} & I_{np} &  \alpha I_{np}
		\end{bmatrix}
		\begin{bmatrix}
			\xi_x \\
			\xi_v \\
			w
		\end{bmatrix}, \notag
	\end{align}
	where $\mu_1, \mu_2 \in \mathbb{R}^+$. Using the same analysis as previous section, if $\mu_1 \lambda_2 +\mu_2 > 1$ and $\alpha>1$, $V$ is positive definite. The derivative of $V$ is given as
	\begin{align} \label{ini-lyp-der}
		\dot{V}=& 
		\mu_1  \xi_x^T [L \otimes I_p] \xi_v + \mu_2  \xi_x^T \xi_v+  \xi_v^T \xi_v+w^T \xi_v -\kappa \xi_x^T \xi_x  \notag \\
		& -\kappa \xi_x^T \xi_v-\alpha \xi_x^T [L \otimes I_p] \xi_x-\kappa \xi_v^T \xi_x -\kappa \xi_v^T \xi_v  \notag \\
		&-\alpha \xi_v^T [L \otimes I_p] \xi_x -\kappa w^T \xi_x+\kappa w^T [M \otimes I_p] r  \notag \\
		&-\kappa w^T \xi_v+\kappa w^T [M \otimes I_p] v^r + \beta w^T [M \otimes I_p]w    \notag \\
		&+\gamma w^T [M \otimes I_p] \mbox{sgn}( w)+ w^T [M \otimes I_p] a^r  \notag \\
		&+ \xi_x^T [L \otimes I_p] \xi_x-\alpha \xi_x^T [L \otimes I_p] \xi_v+\xi_v^T [L \otimes I_p] \xi_x  \notag \\
		&- \alpha \xi_v^T [L \otimes I_p] \xi_v - \alpha \beta w^T w-\alpha \gamma w^T \mbox{sgn}( w)-\alpha \kappa  w^T r \notag \\
		&-\alpha \kappa  w^T v^r-\alpha w^T a^r-\alpha^2 w^T  [L \otimes I_p] \xi_v.
	\end{align}
	The term $ w^T \{ [M \otimes I_p](\kappa r+\kappa v^r+ a^r+\gamma  \mbox{sgn}( w))  -\alpha(\kappa r+\kappa  v^r+ a^r +\gamma  \mbox{sgn}( w)) \}$ can be analyzed as
	\begin{align} \label{ini-w-bound}
		& w^T \{ [M \otimes I_p](\kappa r+\kappa v^r+  a^r + \gamma  \mbox{sgn}( w) ) -\alpha(\kappa r +\kappa  v^r \notag \\
		&+ a^r + \gamma  \mbox{sgn}( w) ) \}  \notag \\
		\leq &  [(1+\alpha)(\kappa \bar{r}+\kappa \bar{v}^r+a^r) +\gamma \sqrt{n}] \cdot ||w||_2 -\alpha \gamma  ||w||_1 \notag \\
		\leq & [(1+\alpha)(\kappa \bar{r}+\kappa \bar{v}^r+a^r) +\gamma \sqrt{n}] \cdot ||w||_1 -\alpha \gamma  ||w||_1.
	\end{align}
	Combining \eqref{ini-lyp-der} and \eqref{ini-w-bound}, if $\gamma > \frac{\alpha+1}{\alpha-\sqrt{n}}(\kappa \bar{r}+\kappa \bar{v}^r+\bar{a}^r)$ and $\mu_1=2\alpha-1$ and $\mu_2=2 \kappa$, we can get
	\begin{align} \label{ini-lyp-der-final}
		\dot{V} \leq &  \xi_v^T \xi_v+w^T \xi_v -\kappa \xi_x^T \xi_x  +(1-\alpha) \lambda_2(L) \xi_x^T \xi_x -\kappa \xi_v^T \xi_v \notag \\
		& -\kappa w^T \xi_x -\kappa w^T \xi_v+ \beta w^T w - \alpha \lambda_2(L) \xi_v^T \xi_v - \alpha \beta w^T w  \notag \\
		&  -\alpha^2 w^T  [L \otimes I_p] \xi_v \notag \\
		=& \begin{bmatrix}
			\xi_x^T &
			\xi_v^T &
			w^T
		\end{bmatrix}
		\begin{bmatrix}
			P_{11} &  P_{12} & P_{13}  \\
			P_{21} &  P_{22} & P_{23}  \\
			P_{31} &  P_{32} & P_{33}
		\end{bmatrix}
		\begin{bmatrix}
			\xi_x  \\
			\xi_v  \\
			w
		\end{bmatrix},
	\end{align}
	where $P_{11}=[-\kappa+(1-\alpha) \lambda_2(L) ] I_{np}$, $P_{22}=[1-\kappa-\alpha \lambda_2(L)]I_{np}$, $P_{12}=P_{21}=0_{np}$ $P_{13}=P_{31}=-\frac{1}{2}\kappa I_{np}$, $P_{23}=P_{32}=\frac{1}{2}[(1-\kappa)I_{np}-\alpha^2 L] $, and $P_{33}=\beta(1-\alpha)I_{np}$ and we have used the fact that $M-I_n$ is negative semi-definite and Lemma \ref{eigen} for the inequality.
	
	If the control gains $\kappa$, $\alpha$, and $\beta$ satisfy the constraints mentioned in Theorem \ref{thm:ini-free-DAC}, the matrix $P$ and thus  $\dot{V}$ are negative definite. Integrating both side of \eqref{ini-lyp-der-final}, it is concluded that $\xi_x$, $\xi_v$ and $w \in \mathbb{L}_2$. Since $V$ is positive definite and $\dot{V}$ is negative definite, $\xi_x$, $\xi_v$ and $w \in \mathbb{L}_\infty$.
	Since $\dot{\xi}_1$, $\dot{\xi}_2$ and $\dot{w}$ are bounded, it is concluded that $\xi_x$, $\xi_v$, $w \to 0$ asymptotically. Therefore it is proved that ${x_i-\frac{1}{n} \sum_{j=1}^n x_j \to 0}$, ${v_i-\frac{1}{n} \sum_{j=1}^n v_j \to 0}$ and $ {w_i \to 0}$ as $t \to \infty$.
	
	In second part, we prove that ${\sum_{i=1}^n x_i \to \sum_{i=1}^n r_i}$ and ${\sum_{i=1}^n v_i \to \sum_{i=1}^n v_i^r}$ asymptotically. Define the variables ${S_1=\sum_{i=1}^n x_i-\sum_{i=1}^n r_i}$ and ${S_2=\sum_{i=1}^n v_i-\sum_{i=1}^n v_i^r}$, we can get from \eqref{ref-dynamic} and \eqref{ini-free-control} that
	\begin{align} \label{initial-value-constraint}
		\dot{S}_1=&S_2, \notag \\
		\dot{S}_2=& -\kappa S_2-\kappa S_1+\beta \sum_{i=1}^n  w_i  +\gamma \sum_{i=1}^n  \mbox{sgn}(w_i).
	\end{align}
	We then use input-to-state stability to analyze the system \eqref{initial-value-constraint} by treating the term $\beta \sum_{i=1}^n  w_i  +\gamma \sum_{i=1}^n  \mbox{sgn}(w_i)$ as the input and $S_1$ and $S_2$ as the states.
	If $\kappa>1$, the matrix $\begin{bmatrix}
	0_p & I_p \\
	-\kappa I_p & -\kappa I_p
	\end{bmatrix}$ is Hurwitz. The system \eqref{initial-value-constraint} with zero input is exponentially stable and hence input-to-state stable. Since $w_i \to 0$ as $t \to \infty$ for each agent, it follows that $S_1 \to 0$ and $S_2 \to 0$, which implies that $\sum_{i=1}^n x_i \to \sum_{i=1}^n r_i$ and $\sum_{i=1}^n v_i \to \sum_{i=1}^n v_i^r$, respectively. Employing the result of the first part, it is concluded that $x_i \to \frac{1}{n} \sum_{j=1}^n r_j$ and ${v_i \to \frac{1}{n} \sum_{j=1}^n v_j^r}$ and the distributed average tracking is achieved.
	\endproof
	
	\begin{remark} Note that the algorithm in \cite{feidoubleintegrator} relies on correct initialization and relative velocity measurements. Using the algorithm defined by \eqref{ini-free-fil-in}-\eqref{ini-free-control}, there is no need for accurate initialization and absolute instead of relative velocity measurements are sufficient. The corresponding trade-off is that this algorithm relies on the assumption of bounded input signals, input velocity and input acceleration.
	\end{remark}
	
	\begin{figure}
		\centering
		\includegraphics[width=0.5\textwidth]{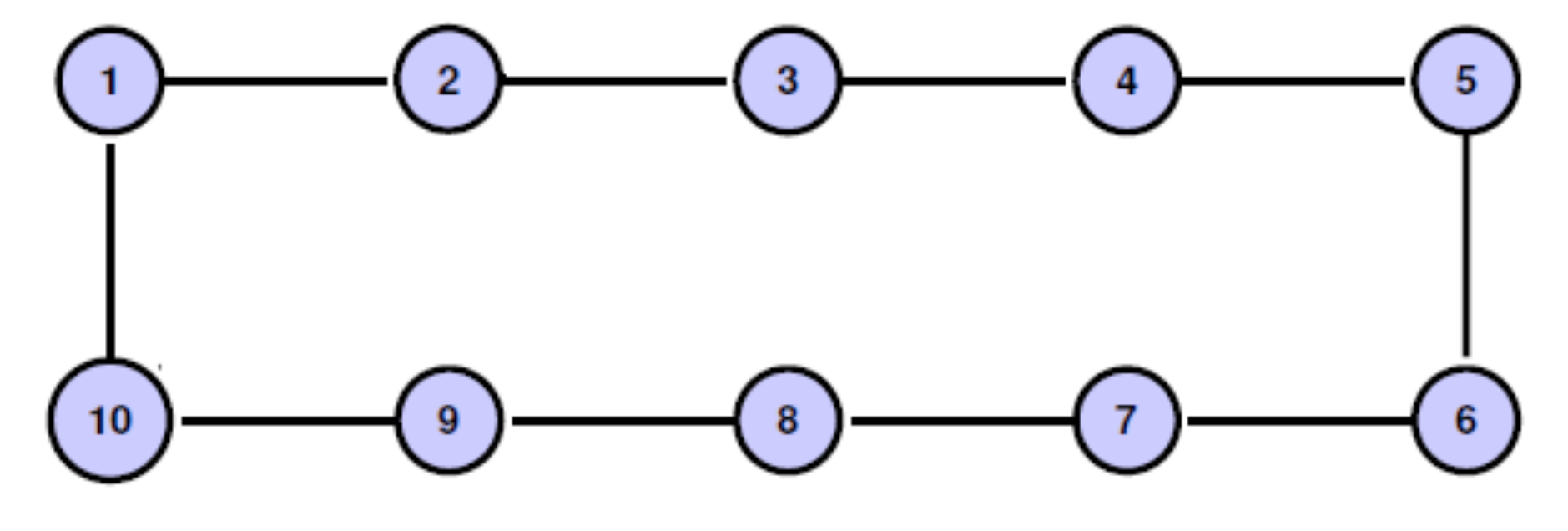}
		\caption{The network topology is used to characterize
			the interaction among ten agents.}
		\label{top}
	\end{figure}
	
	\begin{remark}
		To implement the algorithm \eqref{vel-free-fil-in}-\eqref{vel-free-control}, each agent needs its own and neighbors' positions and filter outputs as well as its own input acceleration.
		Hence communication is needed between neighbors, while there is no need for either absolute or relative velocity measurements. This algorithm needs correct initialization but can deal with more general input signals. To implement the algorithm \eqref{ini-free-fil-in}-\eqref{ini-free-control}, each agent needs its own position, velocity, and relative position between itself and neighbors, and its own input signal, input velocity and input acceleration.
		In this case, communication is not required and the algorithm can rely on only local sensing, where there is no need for relative velocity measurements among neighbors.
		The second one does not use correct initialization but requires the input signal, its velocity and acceleration all to be bounded.
		As a result, both algorithms have their values depending on the application scenarios. Reduced requirement on velocity measurements is achieved in both algorithms to reduce costs.
	\end{remark}

	\section{Simulation}
	
\begin{figure}
		\centering
		\includegraphics[width=0.5\textwidth]{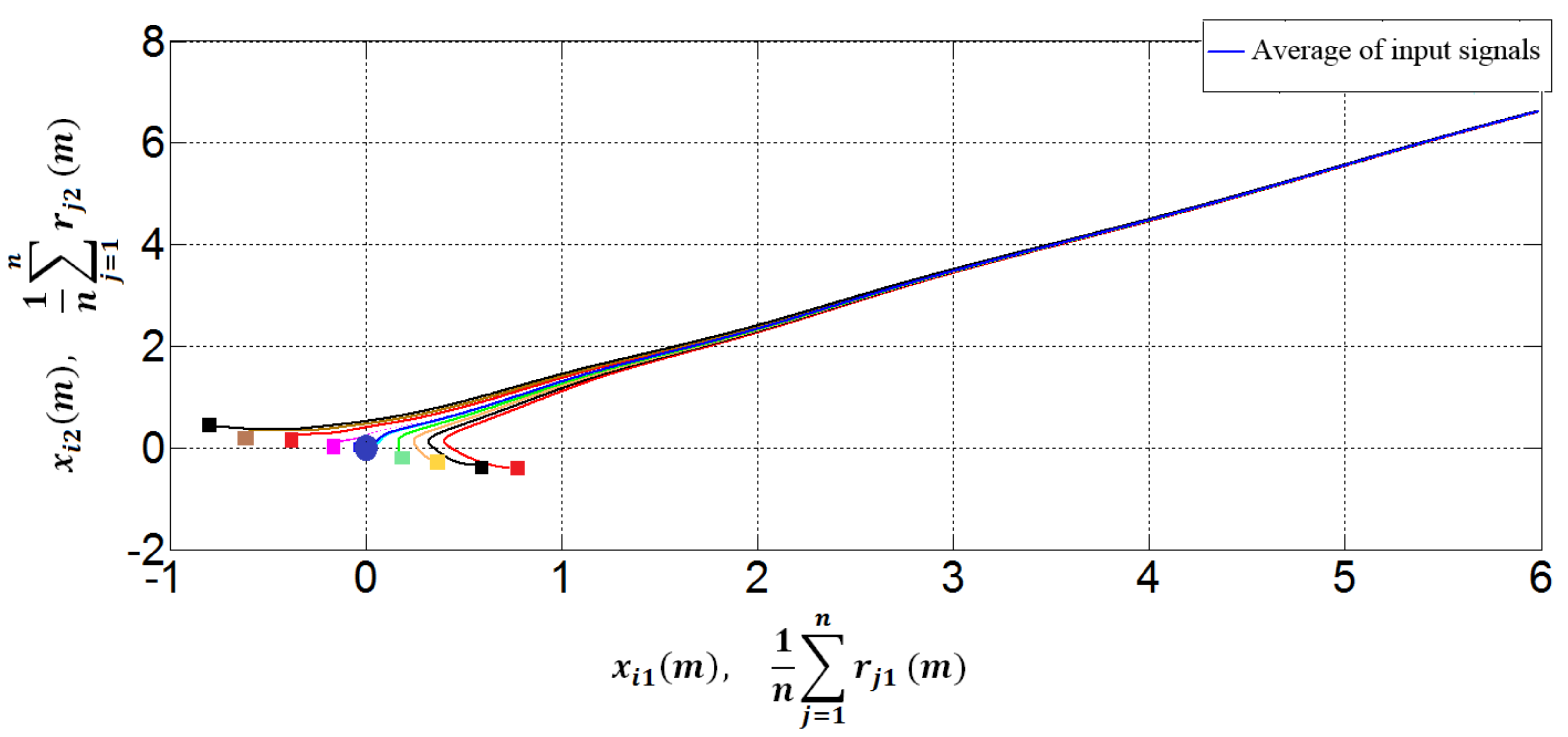}
		\caption{The positions of 10 agents and the average of the input signals, where the states are initialized correctly. The initial average of input signals is represented as a dot while the initial positions of the agents are represented as squares.}
		\label{x2d}
	\end{figure}

\begin{figure}
		\centering
		\includegraphics[width=0.5\textwidth]{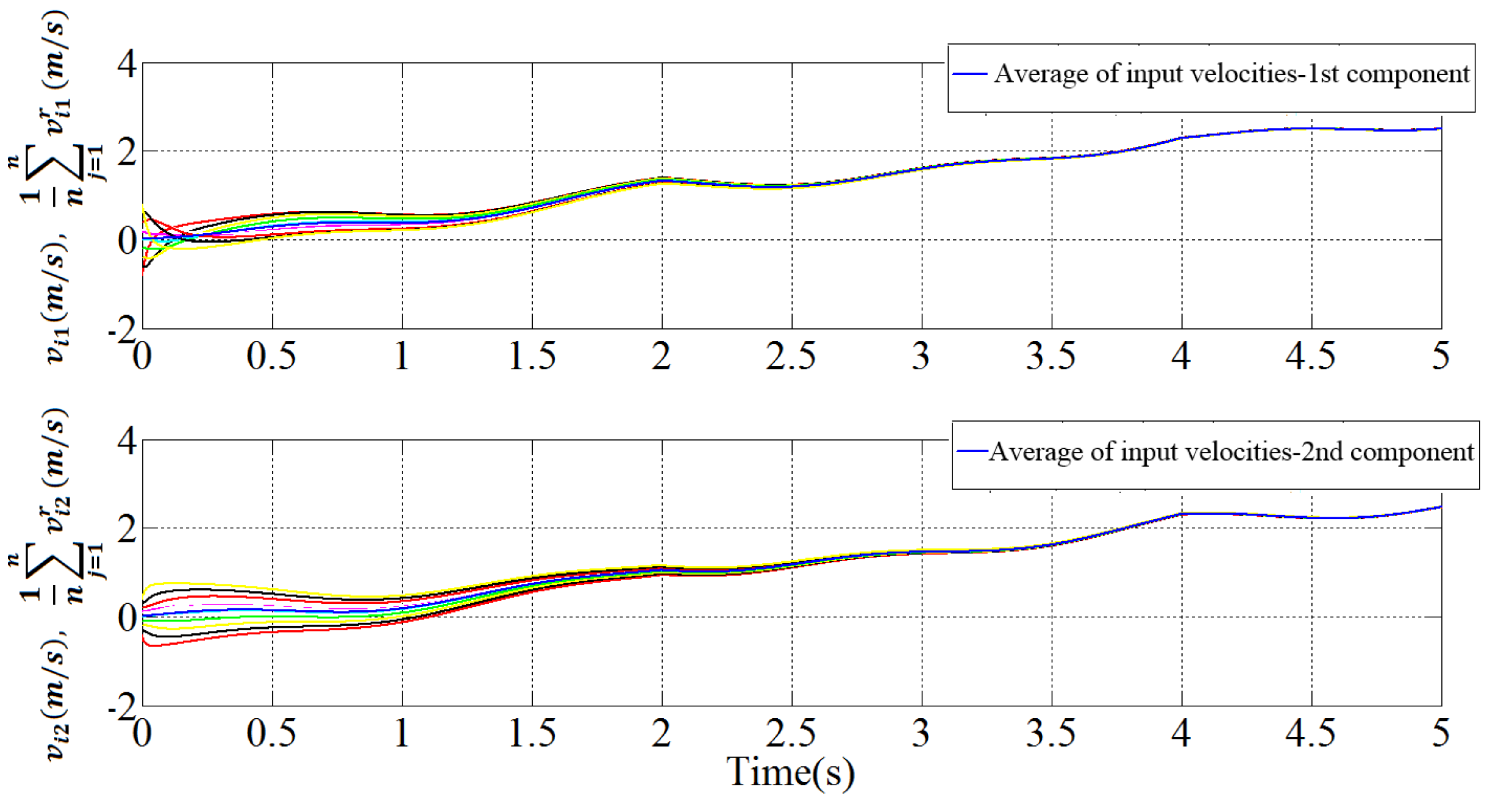}
		\caption{The velocities of 10 agents and the average of input signals, where the agents' velocity are initialized correctly.}
		\label{v2d}
	\end{figure}

	In this section, numerical simulation results are given to illustrate the effectiveness of the theoretical results obtained in Sections \ref{corr-ini} and \ref{abs-corr-ini}.
	It is assumed that there are ten agents $(n = 10)$, where the network topology is described by Fig. \ref{top}.
	In the first case, the input acceleration for agent $i$ is given by $[0.1 (\sin(5t) + \text{mod} (t,2)) \times i,0.1 (\cos(5t) + \text{mod}(t,2)) \times i]^T$, and the initial position and velocity of the agents are chosen as $x(0)=[0.1 \times b,-0.2 \times b]^T$ and ${v(0)=[0.2 \times b,0.1 \times b]^T}$, where $b=[-4,-3,-2,-1,0,1,2,3,4,5]^T$.
	We let $r_i(0)=x_i(0)$ and $v_i^r(0)=v_i(0)$, $i=1,\ldots,10$.
	We denote the $j$th component of $x_i$ as $x_{ij}$. Similar notations are used for $v_i$, $r_i$, and $v_i^r$.
	The control parameters for all agent are chosen as $\alpha=20$, $\gamma=5$ and $\beta=400$.
	We simulate the algorithm defined by \eqref{vel-free-fil-in}-\eqref{vel-free-control}. Fig. \ref{x2d} shows the positions of the agents and the average of the input signals. Clearly, all agents have tracked the average of the input signals in the absence of velocity measurements. Fig. \ref{v2d} shows the velocities of the agents and the average of the input velocities. We see that the distributed average tracking is achieved for the agents' velocities too.

\begin{figure}
		\centering
		\includegraphics[width=0.5\textwidth]{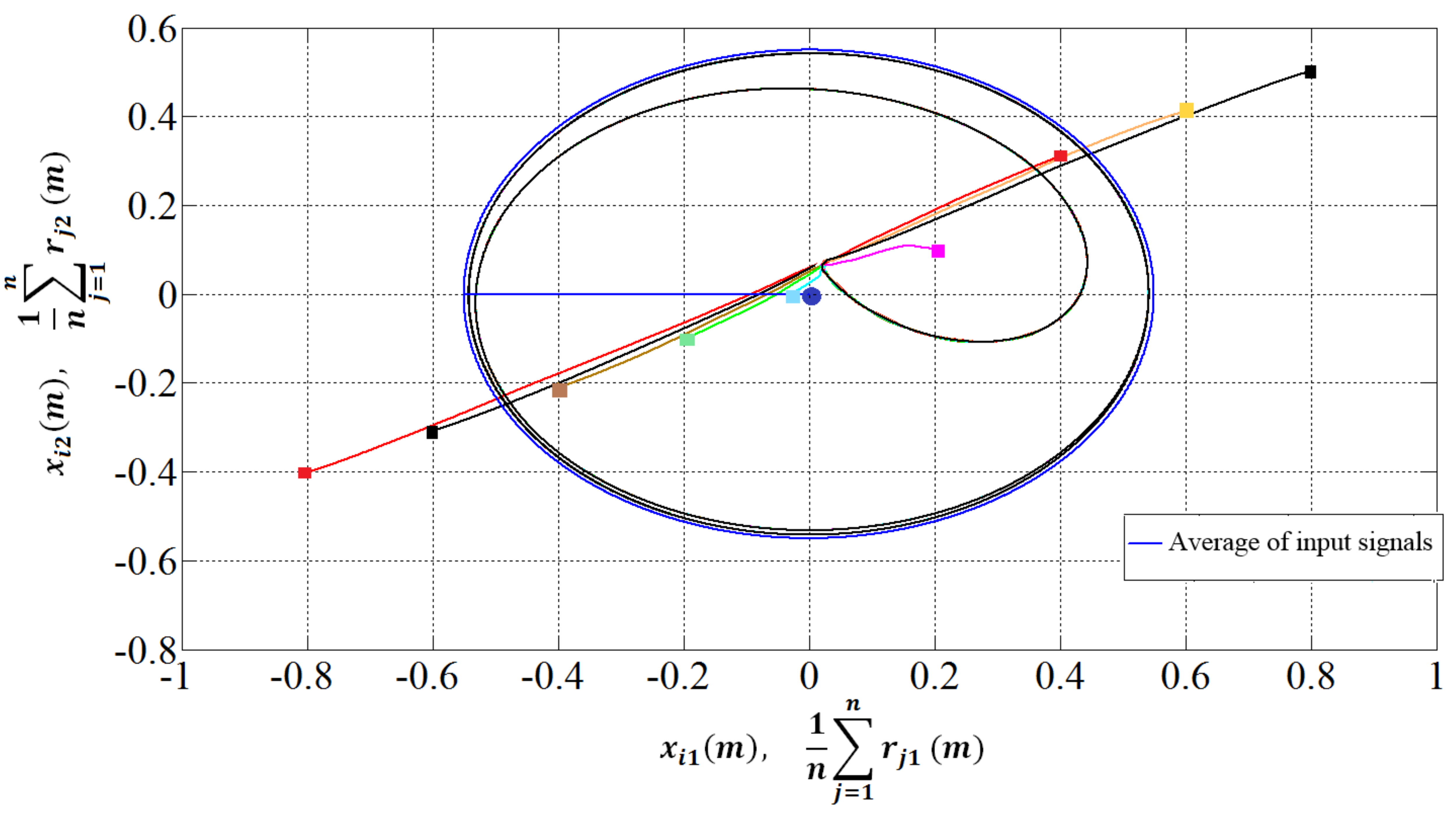}
		\caption{The positions of 10 agents and the average of input signals in the absence of correct initialization. The initial average of input signals is represented as a dot while the initial positions of the agents are represented as squares.}
		\label{x2d2}
	\end{figure}

\begin{figure}
		\centering
		\includegraphics[width=0.5\textwidth]{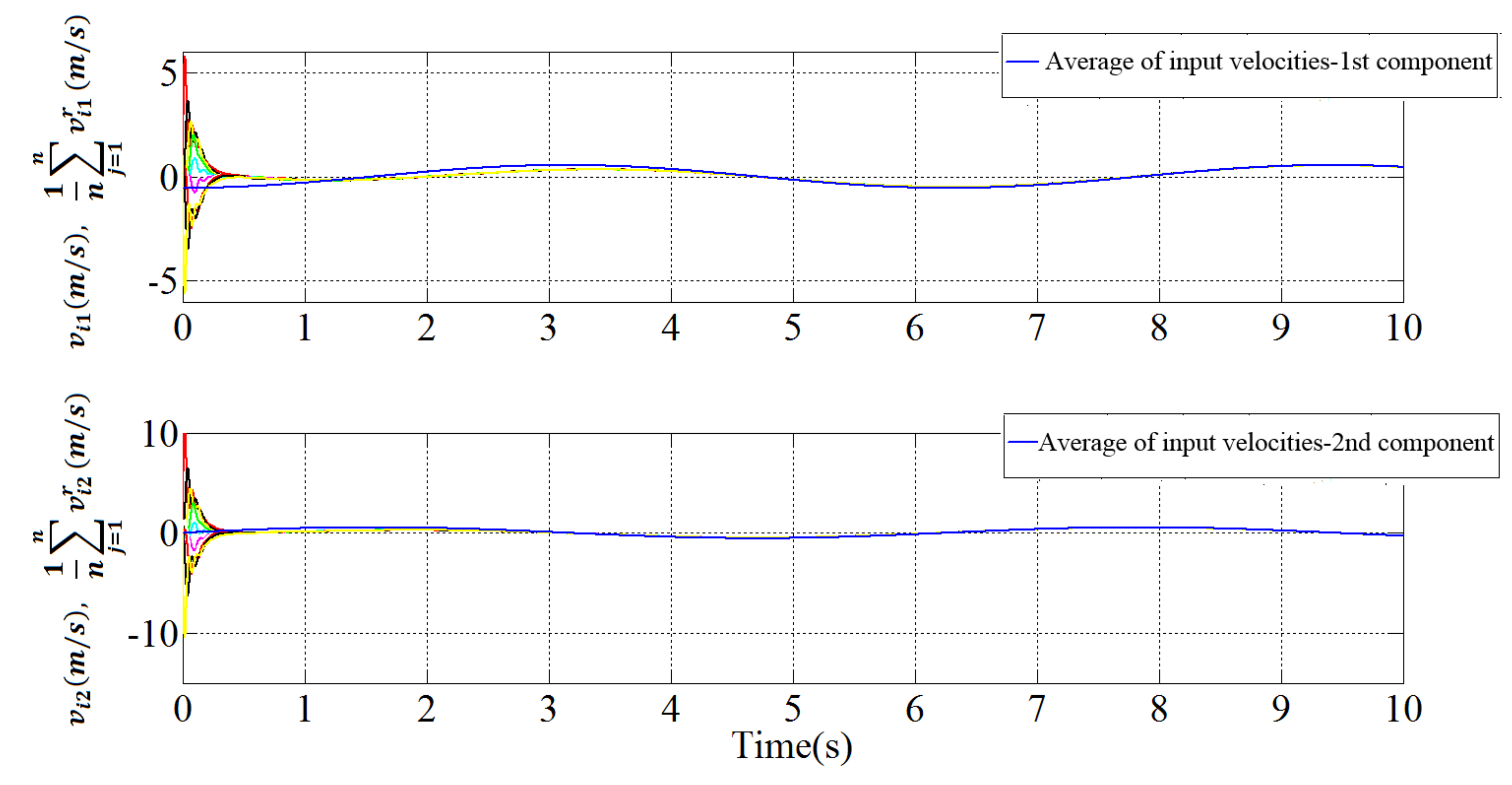}
		\caption{The velocity of 10 agents and the average of input signals in the absence of correct initialization.}
		\label{v2d2}
	\end{figure}

	In the second case, we simulate distributed average tracking when the agents' positions and velocities are not initialized correctly. The input acceleration for each agent is described by ${[0.1 \sin(t) \times i,0.1 \cos(t) \times i]^T}$. We simulate the algorithm defined by \eqref{ini-free-fil-in}-\eqref{ini-free-control}. The initial values are set as ${x(0)=[0.1 \times b,-0.2 \times b]^T}$ and ${v(0)=[0.2 \times b,0.1 \times b]^T}$, $r_i(0)=[0,0]^T$ and ${v_i^r(0)=[-0.1 \times i,-0.1 \times i]^T}$, $i=1,\ldots,10$, and the control parameters are chosen for all agents as $\kappa=2$, ${\alpha=10}$, ${\gamma=50}$ and ${\beta=450}$. Fig. \ref{x2d2} and Fig. \ref{v2d2} show that the distributed average tracking is achieved for both agents' positions and velocities in the absence of neighbors' velocity information and accurate initialization.
	
	\section{CONCLUSIONS}
	
	
	In this paper, distributed average tracking for a group of double-integrator agents has been studied. First a distributed discontinuous controller combined
	with a distributed filter was proposed in the absence of velocity information to deal with input signals with bounded acceleration deviations, where the agents' position and velocities must be initialized correctly.
	Here only the relative position and filter output information were used in the control design besides each agent's own position and input acceleration.
	The algorithm was then modified to remove the requirement of communication and correct initialization under the assumption that the inputs, their velocities and their accelerations are bounded. In both cases, reduced requirement on velocity measurements was achieved and hence the cost for velocity measurements was reduced.
	
	\bibliographystyle{IEEEtran}
	\bibliography{refs}
	
	

\end{document}